\documentclass[12pt]{article}

\usepackage[english]{babel}%
\usepackage[numbers]{natbib}%
\usepackage{enumitem}%
\usepackage{varioref}%
\usepackage{array}%
\usepackage{proof}%
\usepackage{tikz}%
\usepackage{amsmath}%
\usepackage{amsthm}%
\usepackage{amsfonts,amssymb}%
\usepackage[pdftex,pdfusetitle,hidelinks]{hyperref}%

\newcommand{\textcite}[2][]{\citet[#1]{#2}}
\newcommand{\autocite}[2][]{\citep[#1]{#2}}

\addto\captionsenglish{%
  \def\sectionname{Section}%
  \def\definitionname{Definition}%
  \def\examplename{Example}%
  \def\conjecturename{Conjecture}%
  \def\theoremname{Theorem}%
  \def\lemmaname{Lemma}%
  \def\corollaryname{Corollary}%
  \def\remarkname{Remark}%
}

\theoremstyle{plain}%
\newtheorem{lemma}{\lemmaname}[section]%
\newtheorem{theorem}{\theoremname}[section]%
\theoremstyle{definition}%
\newtheorem{definition}{\definitionname}[section]%
\theoremstyle{remark}%
\newtheorem{remark}{\remarkname}[section]%
\newtheorem{example}{\examplename}[section]%

\labelformat{section}{\sectionname~#1}%
\labelformat{subsection}{\sectionname~#1}%
\labelformat{subsubsection}{\sectionname~#1}%
\labelformat{chapter}{\chaptername~#1}%
\labelformat{definition}{\definitionname~#1}%
\labelformat{example}{\examplename~#1}%
\labelformat{conjecture}{\conjecturename~#1}%
\labelformat{corollary}{\corollaryname~#1}%
\labelformat{theorem}{\theoremname~#1}%
\labelformat{lemma}{\lemmaname~#1}%
\labelformat{remark}{\remarkname~#1}%
\labelformat{figure}{\figurename~#1}%

\newcommand{\introductionrule}[1]{\ensuremath{#1}\text{I}}%
\newcommand{\eliminationrule}[1]{\ensuremath{#1}\text{E}}%
\newcommand{\discharge}[1]{[#1]}%
\newcommand{\deduction}[1]{\ensuremath{\Pi_{#1}}}%
\newcommand{\linearargument}[2]{\ensuremath{\langle#1, #2\rangle}}%

\title{On Dummett's Pragmatist Justification Procedure}
\author{Herm\'{o}genes Oliveira}
\date{}

\begin{document}

\maketitle

\begin{abstract}
  I show that propositional intuitionistic logic is complete with
  respect to an adaptation of Dummett's pragmatist justification
  procedure.  In particular, given a pragmatist justification of an
  argument, I show how to obtain a natural deduction derivation of the
  conclusion of the argument from, at most, the same assumptions.
\end{abstract}


\section{Introduction}\label{sec:introduction}
Proof-theoretic definitions of validity can be considered as loosely inspired by Wittgenstein's ideas relating \emph{meaning} and \emph{use}.
They attempt to explain of the concept of logical validity in terms of the \emph{deductive use} of the logical constants, as expressed by inference rules.
In this context, Gentzen's investigations into deduction, particularly his calculus of natural deduction, are often used as a starting point for explaining the meaning of the logical constants on the basis of rules governing their use.

In the standard natural deduction calculus \autocite{gentzenULS,prawitzND}, the deductive use of a logical constant is governed by its introduction and elimination rules.
Thus, from a semantic perspective where meaning is explained on the basis of use, the introduction and elimination rules express the canonical manner in which a sentence with a logical constant as main operator is used in a deductive argument: the introduction rules express the canonical use of the sentence as a conclusion, the elimination rules express the canonical use of the sentence as an assumption.
Along these lines, \textcite{dummettPBIL,dummettLBM} proposed that the analysis of the deductive meaning of a logical constant into introduction and elimination rules accounts for two distinct aspects of its use.
Roughly speaking, the introduction rules show how to establish a sentence, or to warrant its assertion; they stand for the \emph{verificationist} aspect.
On the other hand, the elimination rules show what consequences can be extracted from a sentence, or what difference accepting it makes to our practice; they stand for the \emph{pragmatist} aspect.

Accordingly, in a verificationist approach to proof-theoretic semantics, valid arguments are defined on the basis of introductions rules.
The main idea is that an argument is valid if, whenever we can obtain the assumptions in a canonical manner, we can also obtain the conclusion in a canonical manner.
In a pragmatist approach, on the other hand, valid arguments are defined on the basis of elimination rules.
The main idea is that any consequence that can be drawn in a canonical manner from the conclusion can also be drawn in a canonical manner from the assumptions.

\textcite{prawitzIRPT,prawitzAGPTCCILR} and \textcite{dummettPBIL,dummettLBM} conjectured that proof-theoretic approaches to logical semantics would result in an intuitionistic, or constructive, notion of validity.
Due to bias towards the verificationist point of view, proof-theoretic definitions of validity have often been approached via introduction rules \autocite{prawitzIRPT,prawitzTOFVPV,prawitzMAP,prawitzAGPTCCILR,schroederheisterVCPTS}.
Unfortunately, with respect to their adequacy to intuitionistic logic, verificationist proposals ran into some problems \autocite{sandqvistCLWB,sanzetalCSARVPL,piechaetalFCPTS,piechaCPTS,goldfarbDPTJLL}.

However, as suggested by some remarks of \textcite{gentzenULS}, the introduction and elimination rules for a logical constant are harmoniously related, such that one could extract, in some sense, the elimination rules from the introduction rules and vice versa.
The harmony between introduction and elimination rules suggests that the elimination rules could just as well provide a basis upon which to develop a proof-theoretic definition of validity.

The pragmatist proof-theoretic approach to validity, the one based on elimination rules, has received comparatively less attention \autocite{dummettLBM,prawitzPVTM,schroederheisterPTVBER}.
In this paper, I show that propositional intuitionistic logic is complete with respect to an adaptation of Dummett's pragmatist proof-theoretic definition of validity.
I adapt \mbox{Dummett's} definitions lightly in order to avoid problems and objections but otherwise stay as closely as possible to the original framework.

\section{Preliminaries}\label{sec:preliminaries}
This section explains the basic notions and terminology.
It also recalls and restates results and definitions that will be needed later.
\subsection{Terminology and notation}\label{sec:terminology}
I borrow most of the terminology from \textcite{dummettLBM} himself.
A noteworthy difference is the use of the term ``assumption'' in place of Dummett's term ``initial premiss''.
I included this subsection in the interest of self-containment, but it can be safely skipped if you are already familiar with the notation and typographical conventions used in natural deduction.

\paragraph{The language.}
We consider a propositional language with infinitely many propositional variables (atomic sentences) and the propositional logical constants: $\to$ (implication), $\vee$ (disjunction), $\wedge$ (conjunction) and $\bot$ (absurdity).
The complex sentences of the language are formed from atomic sentences by means of composition with the logical constants in the usual way.
Latin letters ($A$, $B$, $C$ etc.) are used to stand for arbitrary sentences of the language and Greek letters ($\Gamma$ and $\Delta$) to stand for finite collections of sentences.
Subscripts are used whenever it is necessary or convenient.
The \emph{degree of a sentence} is the number of logical constants that occur in it.

\paragraph{Arguments and derivations.}
Formally, \emph{arguments} can be considered as trees of sentence occurrences (designated with \deduction{}, possibly with subscripts).
They are constructed from top to bottom, from the leaves to the root, by \emph{inferences}.
These inferences lead from one or more sentences, the \emph{premisses}, to a single sentence, the \emph{consequence}.
In an argument, each premiss of an inference is either a leaf of the tree or the consequence of a previous inference.
Thus, argument trees are formal representations of the process of argumentation, or reasoning, with some leaves acting as \emph{assumptions} and the root acting as the \emph{conclusion} of the argument.
Any occurrence of a sentence in an argument determine, in the obvious way, a subargument with that sentence as conclusion.\footnote{As a limiting case, a single sentence occurrence is an argument with that sentence acting as both assumption and conclusion.}
A \emph{path} in an argument is a sequence of sentence occurrences such that each sentence in the path is an immediate inferential consequence of the previous one.\footnote{Again, as a limiting case, a single sentence occurrence measures an \emph{empty} path from that sentence occurrence to itself.}
Every leaf in an argument is initially an assumption, albeit assumptions can be discharged by inferences.\footnote{An axiom or logical theorem $A$ can be considered the result of an inference from leaf $A$ to conclusion $A$ that discards the leaf occurrence of $A$.}
After an assumption is discharged by an inference, the argument, starting from the consequence of that inference, does not depend any more on the assumption.
The discharge of assumptions are indicated using square brackets with numeric indices used to pinpoint the particular inference discharging the assumption.
Whenever it is clear from context, the numeric indices are left implicit.  I write \linearargument{\Gamma}{A} to denote an argument from assumptions $\Gamma$ (those that remained undischarged throughout the argument) to conclusion $A$ without paying attention to the argumentation process that goes from $\Gamma$ to $A$.
Also, for the sake of simplicity, I often talk about sentences when I actually mean occurrences of sentences in an argument, and similarly with respect to inference rules and the particular inferences resulting from their application.

Propositional intuitionistic logic is characterised by the standard system of natural deduction \autocite{gentzenULS,prawitzND}.
The inference rules for the propositional connectives are symmetrically distributed between introduction (I) and elimination (E) rules.
\begin{displaymath}
  \begin{array}
    {p{0.2\textwidth}p{0.2\textwidth}p{0.2\textwidth}p{0.2\textwidth}}
    \infer[\introductionrule{\to}]{A\to{B}}{
      \infer*{B}{\discharge{A}}}
    & & \infer[\eliminationrule{\to}]{B}{A & A\to{B}} &
  \end{array}
\end{displaymath}
\begin{displaymath}
  \begin{array}
    {p{0.2\textwidth}p{0.2\textwidth}p{0.2\textwidth}p{0.2\textwidth}}
    \infer[\introductionrule{\wedge}]{A\wedge{B}}{A & B} & &
    \infer[\eliminationrule{\wedge}]{A}{A\wedge{B}} &
    \infer[\eliminationrule{\wedge}]{B}{A\wedge{B}}
  \end{array}
\end{displaymath}
\begin{displaymath}
  \begin{array}
    {p{0.2\textwidth}p{0.2\textwidth}p{0.2\textwidth}p{0.2\textwidth}}
    \infer[\introductionrule{\vee}]{A\vee{B}}{A} &
    \infer[\introductionrule{\vee}]{A\vee{B}}{B} &
    \infer[\eliminationrule{\vee}]{C}{A\vee{B} &
      \infer*{C}{\discharge{A}} & \infer*{C}{\discharge{B}} } &
  \end{array}
\end{displaymath}
Negation ($\lnot$) can be defined as usual in terms of implication and absurdity.
The rule for the absurdity logical constant $\bot$ can be considered an elimination rule.
\begin{displaymath}
\infer[\eliminationrule{\bot}]{A}{\bot}
\end{displaymath}
Natural deduction derivations are a particular subclass of arguments in which every inference is in accordance with one of the inference rules above.

\subsection{Elimination rules and related notions}\label{sec:elim-rules}
In an elimination rule for a logical constant, exactly one premiss of the rule is required to have that constant as main logical operator.
This premiss is called the \emph{major premiss}, and all others, if there are any others, are called \emph{minor premisses}.

In an elimination rule, a minor premiss is \emph{vertical} if the same sentence figures as both minor premiss and consequence of the rule, otherwise it is called \emph{horizontal}.
An elimination rule is a \emph{vertical rule} if at least one of its minor premisses is vertical \emph{and} it allows the discharge of assumptions in the
subarguments for its vertical minor premisses.
We can and shall assume that every application of vertical elimination rules do in fact discharge the assumptions as indicated by the rule.\footnote{Otherwise, the application of the elimination rule would be superfluous \autocite[p.~49--50]{prawitzND}.
  This is in line with \textcite[p.~283]{dummettLBM}.}
Elimination rules which are not vertical are called \emph{reductive}.

\subsection{Normal derivations}\label{sec:normalisation}
This subsection recollects some results and fixes some terminology regarding normal derivations in propositional intuitionistic logic.
The results are restated mainly for your convenience and to avoid confusion resulting from conflicting terminology.
For detailed proofs, please resort to Prawitz's monographs \autocite{prawitzND,prawitzIRPT}.
\begin{definition}\label{def:track}
  A \emph{track}\footnote{The notion of ``track'' is obviously a slight adaptation of Prawitz's notion of ``path''.  Since the term ``path'' is already reserved for a different concept, the term ``track'' is used in order to avoid confusion.} in an argument \deduction{} is a sequence of \emph{sentence occurrences} $A_{1}, \cdots{} A_{n}$ such that:
\begin{enumerate}[label=(\roman*)]
\item $A_{1}$ is a leaf in \deduction{} that is not discharged by an application of a vertical elimination rule
\item $A_{i}$, for each $i<n$, is not a horizontal minor premiss of an application of an elimination rule, and either (1) $A_{i}$ is not a major premiss of a vertical rule and $A_{i+1}$ is the sentence immediately below $A$, or (2) $A_{i}$ is a major premiss of a vertical rule and $A_{i+1}$ is the assumption discharged by the respective application of the rule
\item $A_{n}$ is either a horizontal minor premiss or the conclusion of \deduction{}
\end{enumerate}
\end{definition}
As a sequence of occurrences of sentences, a track can also be divided into \emph{segments}; they conflate repeated consecutive occurrences of the same sentence that arise from applications of vertical elimination rules.
If the first sentence in a track is an assumption of the argument, the track is \emph{open} in the argument, otherwise it is \emph{closed}. 

The tracks in an argument can be assigned an order.
The lowest order is assigned to tracks whose last sentence occurrence is the conclusion of the argument; they are called \emph{main tracks}.
The order then increases progressively through horizontal minor premisses.
The major premiss of the horizontal rule belongs to a \emph{parent track} which, of course, can have other \emph{children tracks} sharing an immediately higher order than their parents.
The last sentence occurrence in a track determines a subargument, or subderivation, whose order is the order of the track.
The progeny relationship between tracks can be naturally extended to cover them.
\begin{theorem}\label{thm:normal}
  Let $\tau$ be a track in a normal intuitionistic derivation and let $\sigma_{1}\cdots\sigma_{n}$ be the corresponding sequence of segments in $\tau$.
  Then there is a segment $\sigma_{i}$, called the \emph{base segment} in $\tau$, which separates two (possibly empty) parts of $\tau$, called the \emph{analytic part} and the \emph{synthetic part} of $\tau$, such that:
  \begin{enumerate}[label=(\roman*)]
  \item for each $\sigma_{j}$ in the analytic part, $\sigma_{j}$ is a major premiss of an elimination rule and the sentence occurring in $\sigma_{j+1}$ is a subsentence of the one occurring in $\sigma_{j}$
  \item the base segment $\sigma_{i}$ is a premiss of an introduction rule or of \eliminationrule{\bot}, provided $i\neq{n}$
  \item for each $\sigma_{j}$ in the synthetic part, except the last one, $\sigma_{j}$ is a premiss of an introduction rule and the sentence occurring in $\sigma_{j}$ is a subsentence of the one occurring in $\sigma_{j+1}$
  \end{enumerate}
\end{theorem}
\begin{definition}\label{def:positive-negative}
  The subsentences of a sentence $A$ are classified as \emph{positive} or \emph{negative} as follows:
  \begin{itemize}
  \item $A$ is a positive subsentence of $A$
  \item if $B\wedge{C}$ or $B\vee{C}$ are positive (resp. negative) subsentences of $A$, then $B$ and $C$ are positive (resp. negative) subsentences of $A$ 
  \item if $B\to{C}$ is a positive (resp. negative) subsentence of $A$, then $B$ is a negative (resp. positive) subsentence of $A$ and $C$ is a positive (resp. negative) subsentence of $A$
  \end{itemize}
\end{definition}

\begin{definition}\label{def:components}
  A sentence $A$ is an \emph{assumption component} (resp. \emph{conclusion component}) of an argument \linearargument{\Gamma}{G} when $A$ is a positive (resp. negative) subsentence of some assumption in $\Gamma$, or a negative (resp. positive) subsentence of the conclusion $G$.
\end{definition}
The notions defined above can be naturally extended to cover segments, whereby a segment $\sigma_{j}$ is a subsegment of a segment $\sigma_{i}$  if the sentence occurring in $\sigma_{j}$ is a subsentence of the sentence occurring in $\sigma_{i}$.
\begin{theorem}\label{thm:components}
  Let $\sigma_{1}\cdots{\sigma_{n}}$ be a track in a normal derivation of $G$ from $\Gamma$.
It holds that:
\begin{enumerate}[label=(\roman*)]
\item every segment occurring in the analytic part is an assumption component of \linearargument{\Gamma}{G} and subsegment of $\sigma_{1}$
\item the base segment $\sigma_{i}$ is an assumption component of \linearargument{\Gamma}{G} and a subsegment of $\sigma_{1}$; also, if different from $\bot$, $\sigma_{i}$ is a conclusion component of \linearargument{\Gamma}{G} and a subsegment of $\sigma_{n}$
\item every segment occurring in the synthetic part is a conclusion component of \linearargument{\Gamma}{G} and a subsegment of $\sigma_{n}$
\end{enumerate}
\end{theorem}
\ref{thm:complexity} below provides a complexity measure for normal derivations that decreases on the order of the tracks.
This complexity measure is later employed in \ref{def:validity-canonical} and \ref{def:complementation}.
\begin{definition}\label{def:components-partition}
  A \emph{negative assumption} (resp. \emph{conclusion}) \emph{component} of an argument \linearargument{\Gamma}{G} is a negative subsentence of $G$ (resp. some sentence in $\Gamma$).
Analogously, a \emph{positive conclusion} (resp. \emph{assumption}) \emph{component} of an argument \linearargument{\Gamma}{G} is a positive subsentence of $G$ (resp. some sentence in $\Gamma$).
\end{definition}
The negative and positive partitions introduced in \ref{def:components-partition} make out the components in \ref{def:components}.
The assumption and conclusion components afford a coarse summary of the sentences occurring in the elimination and introduction parts of tracks in a normal derivation on the basis of its assumptions and conclusion.
Their partition into negative and positive constituents yields a somewhat finer distinction.
\begin{definition}\label{def:distinctive-components}
  A \emph{distinctive assumption} (resp. \emph{conclusion}) \emph{component} of an argument is either a negative assumption (resp. conclusion) component or an assumption (resp. the conclusion).
\end{definition}
The distinctive assumption components correspond intuitively to first sentences (leaves) of tracks.
They are either assumptions or, possibly, leaves discharged by applications of \introductionrule{\to} and, consequently, negative subsentences of the conclusion.
The negative assumption components represent, in a sense, a \emph{potential} discharge.
Whether an actual leaf was discharged would depend on the particular derivation.

The distinctive conclusion components correspond intuitively to last sentences (roots) of tracks.
They include the conclusion of main tracks and, possibly, conclusions of their descendant tracks.
More precisely, provided a parent track is open, the conclusion of a child track is a negative subsentence of its assumption.
\begin{definition}[\textcite{dershowitzPTMO}]\label{def:degree-collection}
  The \emph{degree of a collection of sentences} $\Delta_{j}$ is lower than the degree of a collection of sentences $\Delta_{i}$ ($\Delta_{j}<\Delta_{i}$) if, and only if, $\Delta_{j}$ results from $\Delta_{i}$ by replacing one or more sentences with a finite collection of sentences of lower degree. 
\end{definition}
\begin{lemma}\label{lem:conclusion-component}
  In a normal derivation, the degree of the collection of distinctive conclusion components never increases on the order of subderivations.
  In particular, if a parent track is closed, then the collection of distinctive conclusion components of the child subderivation has lower degree.
\end{lemma}
\begin{proof}
  Consider derivations \deduction{i} and \deduction{j}, where \deduction{j} is a child of \deduction{i}.
  I show that any distinctive conclusion component of \deduction{j} is either itself a distinctive conclusion component of \deduction{i}, or is replaced in \deduction{i} with a distinctive conclusion component of higher degree.
  By \ref{def:distinctive-components}, a distinctive conclusion component of \deduction{j} is either the conclusion of \deduction{j} or a negative conclusion component of \deduction{j}.
  First, consider the negative conclusion components of \deduction{j}.
  By \ref{def:components-partition}, they also belong to the distinctive conclusion components of \deduction{i}, unless the corresponding assumption was discharged.
In that case, by \ref{thm:normal},  either (\introductionrule{\to}) the conclusion of \deduction{i} has higher degree, or (\eliminationrule{\vee}) it is a negative subsentence of some assumption in \deduction{i} and, consequently, a negative conclusion component of \deduction{i}.
Now consider the conclusion of \deduction{j}.
Take a parent track $\tau$ in \deduction{i}.
Suppose that $\tau$ is open and $A$ is its assumption.
By \ref{thm:components}, it is a negative subsentence of $A$ and, by \ref{def:distinctive-components}, also a distinctive conclusion component of \deduction{i}.
Finally, suppose that $\tau$ is closed and its leaf $A$ was discharged.
Then, by \ref{thm:normal}, the conclusion of \deduction{i} has higher degree than the conclusion of \deduction{j}.
\end{proof}
\begin{definition}\label{def:complexity}
Let $\Delta_{\Gamma}$ (resp. $\Delta_{G}$) stand for the collection of distinctive assumption (resp. conclusion) components of a normal derivation from $\Gamma$ to $G$.
Consider sequences of subderivations $\linearargument{\Gamma_{0}}{G_{0}}\cdots\linearargument{\Gamma_{n}}{G_{n}}$ where each element stands for a child subderivation of the previous one (\linearargument{\Gamma_{0}}{G_{0}} being the main derivation).
The derivations in the sequence are measured for their \emph{complexity} as follows:
  \begin{itemize}
  \item if $\Delta_{G_{j}}<\Delta_{G_{i}}$ then $\linearargument{\Gamma_{j}}{G_{j}}\prec\linearargument{\Gamma_{i}}{G_{i}}$
  \item if $\Delta_{G_{j}}=\Delta_{G_{i}}$ then
    \begin{itemize}
    \item if $\Delta_{\Gamma_{j}}<\Delta_{\Gamma_{i}}$ then $\linearargument{\Gamma_{j}}{G_{j}}\prec\linearargument{\Gamma_{i}}{G_{i}}$
    \item if $\Delta_{\Gamma_{j}}\geq\Delta_{\Gamma_{i}}$ then $\linearargument{\Gamma_{j}}{G_{j}}\succeq\linearargument{\Gamma_{i}}{G_{i}}$
    \end{itemize}
  \end{itemize}
\end{definition}
\begin{theorem}\label{thm:complexity}
  For any derivable argument, there is a normal derivation where the complexity decreases on the order of subderivations.
\end{theorem}
\begin{proof}
  I describe a method to shorten normal derivations by removing redundancies (loops).
  Consider a normal derivation where a subderivation \deduction{j} has equal or higher complexity than its parent subderivation \deduction{i}.
By \ref{lem:conclusion-component}, a parent track is open.
Let $A_{i}$ be its assumption.
Now, there is also an open track in \deduction{j} where another occurrence $A_{j}$ of the same sentence is an analytic assumption, because otherwise \deduction{j} would already have lower complexity than \deduction{i}.
Let \deduction{k} be the descendant subderivation for the minor premiss of $A_{j}$.
Replace \deduction{j} with \deduction{k}.
Because \deduction{j} has at least the same complexity as \deduction{i}, any negative assumption components of \deduction{j} is also a negative assumption component of \deduction{i}.
Therefore, any discharge of assumptions of \deduction{k} can be transferred from \deduction{j} to \deduction{i}.
This shortening process can be iterated until the resulting child of \deduction{i} has lower complexity.
\end{proof}
\section{Proof-theoretic validity}\label{sec:validity}
The proof-theoretic justification procedures proposed by \textcite{dummettLBM} consist of definitions of validity for arguments based on canonical inference rules for the logical constants.
These inference rules fix the meaning of the logical constants by expressing their canonical deductive use.  They are, in Dummett's terminology, self-justifying.

In contrast with common practice \autocite{prawitzIRPT,prawitzTFGPT,schroederheisterVCPTS,schroederheisterPTVBER}, Dummett's definitions are not based on semantic clauses for particular logical constants.
Instead, he assumes that self-justifying rules are given.
These self-justifying rules are introduction rules in the context of the verificationist procedure, and elimination rules in the context of the pragmatist procedure.
In both procedures, the definitions are stated irrespective of the particular constants or rules provided.
\mbox{Dummett's} approach is thus more general and could, in principle, be applied to any logic.\footnote{For technical reasons, the approach would certainly be restricted at the outset to logics with nice proof-theoretic properties (normalisation, subformula property and etc.), however there is no intrinsic technical limitation that confines it to intuitionistic logic or some specific formulation thereof.}
For our particular case, the self-justifying rules are the standard elimination rules of propositional intuitionistic logic: \eliminationrule{\wedge}, \eliminationrule{\to}, \eliminationrule{\vee} and \eliminationrule{\bot}.

From a pragmatist perspective, a canonical argument start from open (usually complex) assumptions and, through applications of elimination rules, arrive at an atomic conclusion.
Dummett also admits \emph{basic rules} (or \emph{boundary rules}) to determine deducibility among atomic sentences.
In canonical arguments, these basic rules can be applied to an atomic conclusion in order to obtain further atomic consequences.

However, since our main concern is with logical validity, I leave basic rules out of the picture and adapt Dummett's definitions accordingly, for the sake of simplicity.\footnote{\textcite[p.~273]{dummettLBM} explicitly stated the irrelevance of basic rules, or boundary rules, as he called them, to logical validity.
  Furthermore, by reflection on the pragmatist definition of validity (\ref{def:validity}), it is easy to see that any basic rule in the complementation can be transferred to the valid canonical argument required, thus making no difference to which complex arguments are actually validated.}
I also adapt the definitions to the propositional case.
\begin{definition}\label{def:principal}
  A sentence occurrence in an argument is \emph{principal} if every sentence (inclusive) in the path down to the conclusion (exclusive) is a major premiss of an elimination rule.\footnote{As a limiting case, in an argument consisting of a single occurrence of a sentence $A$, acting both as assumption and conclusion, $A$ \emph{is} principal, since the empty path from $A$ to $A$ satisfies the definition.}
\end{definition}
\begin{definition}\label{def:proper}
  An argument is \emph{proper} if at least one of its assumptions is principal.
\end{definition}

The concept of proper argument is an essential component in the pragmatist notion of validity because proper arguments are built from the principal assumption by application of elimination rules.
Arguments that do not follow this pattern are \emph{improper}.
The path from the principal assumption to the conclusion is the \emph{principal path}.

The notion of canonical argument to be introduced later (\ref{def:canonical}) is based primarily on the notion of a proper argument.
Even in proper arguments, however, the subargument for minor premisses of elimination rules may depend on auxiliary assumptions that arrive at the conclusion through improper means, i.\,{}e., through a path that is not solely composed of major premisses of elimination.
These kind of improper subarguments for minor premisses of elimination are called \emph{critical subarguments} (\ref{def:critical}).

\subsection{Canonical arguments and critical subarguments}\label{sec:canonicity}
The following definitions are adapted from \mbox{Dummett's} original definitions as explained in \ref{sec:prawitz-objection}.
The core ideas, however, are preserved.
\begin{definition}\label{def:placid}
  A sentence occurrence is \emph{placid} if no sentence down the path
  to the conclusion is a horizontal minor premiss.
\end{definition}
\begin{definition}\label{def:canonical}
  A \emph{canonical} argument has the following properties:
  \begin{enumerate}[label=(\roman*)]
  \item \label{cond-proper} it is proper;
  \item \label{cond-placid} the subargument for any placid vertical minor
    premiss of an elimination rule is proper.
  \end{enumerate}
\end{definition}
\begin{definition}\label{def:critical}
  A \emph{critical subargument} of a canonical argument is a non-canonical subargument whose conclusion is a horizontal minor premiss of an elimination rule.
\end{definition}

In subarguments for minor premisses, the notion of canonicity deals with vertical and horizontal minor premisses differently.
The subarguments for \emph{vertical} minor premisses are considered independent auxiliary arguments and are thus required to be proper themselves.
Ideally, the subarguments for \emph{horizontal} minor premisses would also be proper (and canonical).
However, \emph{in general}, it is not possible to place any restrictions on the form of the subarguments for horizontal minor premisses: when not already canonical, those subarguments are critical, which means that the validity of the whole canonical argument would depend on their validity (\ref{def:validity-canonical}).
\begin{remark}\label{rmk:critical}
The inference steps in canonical arguments consist primarily of applications of elimination rules: there is a principal path (which is composed of eliminations) and subarguments for minor premisses which are themselves either proper (again with a principal path composed of eliminations) or critical.
Thus, a canonical argument has the general form
\begin{displaymath}
  \infer{\cdots}{\infer{\cdots}{\infer{\cdots}{\cdots} &
      \infer{\cdots}{\infer{\cdots}{\cdots{} & \bigtriangledown{}{}} &
        \bigtriangledown{}{}} &
      \infer{\cdots{}}{\cdots{} & \bigtriangledown{}{}}} &
    \bigtriangledown{}{}}
\end{displaymath}
where the inference steps are applications of elimination rules except for the critical subarguments (represented with ``$\bigtriangledown{}{}$'' above), because the definitions impose no restrictions on their inference steps.
If we were to ignore the critical subarguments, what remained could be called the \emph{proper fraction} of the canonical argument and the corresponding sentence occurrences be called \emph{proper occurrences}.
In the proper fraction, in addition to the principal assumption, all other assumptions are principal assumptions of proper subarguments.
In the context of canonical arguments, they are called, collectively, \emph{proper leaves}, or \emph{proper assumptions}, when undischarged throughout the argument. 
\end{remark}
\begin{lemma}\label{lem:subsentence}
  The conclusion of a canonical argument is always a subsentence of some assumption, provided there is no proper occurrence of $\bot$.
\end{lemma}
\begin{proof}
  Let \deduction{} be a canonical argument with no proper occurrence of $\bot$.
By \ref{def:canonical}, \deduction{} is proper and, by \ref{def:proper} and \ref{def:principal}, it has a principal path of major premisses of applications of elimination rules from an assumption to the conclusion.
In the principal path, the consequences of applications of \eliminationrule{\wedge} and \eliminationrule{\to} are subsentences of their respective major premisses.  The interesting cases are applications of vertical rules (\eliminationrule{\vee}) since they figure a consequence which is not required to be a subsentence of the major premiss.
By \ref{def:canonical}, the subarguments for minor premisses of vertical rules are proper.  By \ref{def:proper} and \ref{def:principal}, each vertical subargument has a path of eliminations from a proper leaf to the conclusion of the subargument.
Now, if the proper leaf of a subargument for a vertical minor premiss was discharged by the corresponding application of \eliminationrule{\vee}, then the conclusion is, by induction hypothesis, a subsentence of its major premiss and hence a subsentence of the principal assumption.
Otherwise, if the proper leaf was not discharged, then it is actually a proper assumption of the canonical argument and the conclusion is, by induction hypothesis, a subsentence of this assumption.
\end{proof}
\begin{remark}\label{rmk:track}
It would perhaps be useful to make a parallel between the definitions above and concepts familiar from normalisation for intuitionistic natural deduction (\ref{sec:normalisation}).
For instance, notice that clause \ref{cond-placid} of \ref{def:canonical} ensures that the segments in main tracks of canonical arguments are all major premisses of applications of elimination rules (except the last one).
In main tracks of canonical arguments, the major premisses of vertical rules are followed by the corresponding discharged assumption, which (if not the last segment in the track) is also a major premiss of an elimination rule.
As a result, the first sentence in a main track is a proper assumption and the last sentence (the conclusion of the canonical argument) is a subsentence of this assumption, provided $\bot$ does not occur in the main track.
The first sentence in a main track, however, can be distinct from the principal assumption of the canonical argument, because assumptions discharged by vertical rules need not be principal in the proper subargument for a minor premiss.
The tracks to whom the principal assumption belongs may be called \emph{principal tracks}.
By \ref{thm:normal}, in normal derivations with empty synthetic parts, the main tracks are all principal tracks and the derivations are, therefore, canonical arguments.
\end{remark}
\begin{definition}\label{def:validity-canonical}
  A canonical argument is \emph{valid} if all its critical subarguments are valid and of lower complexity.
\end{definition}
\begin{definition}\label{def:complementation}
  A \emph{complementation} of an argument \linearargument{\Gamma}{G} is the result of replacing $G$ by a valid canonical argument with the following properties:\footnote{Dummett's original definition has a special clause for when the conclusion $G$ is an atomic sentence.
His clause is subsumed here by canonical arguments consisting of a single occurrence of $G$ (\ref{def:principal}).}
  \begin{enumerate}[label=(\roman*)]
  \item \label{cond-principal} it has $G$ as principal assumption
  \item \label{cond-atomic} it has an atomic conclusion
  \item \label{cond-complexity} it has at most the same complexity
  \end{enumerate}
\end{definition}
\begin{definition}\label{def:validity}
  An argument is \emph{valid} if there is an effective method to transform any complementation of it into a valid canonical argument for the same conclusion from, at most, the same assumptions.
\end{definition}

\ref{def:validity-canonical}, \ref{def:complementation} and \ref{def:validity} should always be considered together since they define notions in terms of each other.  For instance, the notion of valid canonical argument in \ref{def:validity-canonical} is defined in terms of the notion of valid argument which is itself only defined in \ref{def:validity}.

The process of complementation consists basically in the application of elimination rules to the conclusion of the argument until we reach an atomic sentence.
During complementation, the application of elimination rules figuring minor premisses can introduce auxiliary assumptions.
Thus, the valid canonical argument required by \ref{def:validity} can depend on additional assumptions introduced by complementation.
\begin{theorem}[Completeness of Intuitionistic Logic]\label{thm:completeness}
  If an argument \linearargument{\Gamma}{G} is valid, then there is a natural deduction derivation of $G$ from $\Gamma$ in intuitionistic logic.
\end{theorem}
\begin{proof}
Suppose \linearargument{\Gamma}{G} is valid.
By \ref{def:validity}, for any complementation \deduction{c}, we have a valid canonical argument \deduction{v} with, at most, the same assumptions and the
same conclusion.
By \ref{def:complementation}, the complementations are obtained by replacing $G$ with a valid canonical argument that has $G$ as principal assumption.
By \ref{def:proper}, there is a principal path in \deduction{c} from $G$ to the atomic conclusion $C$ which consists solely of applications of elimination rules.
Furthermore, there can be auxiliary assumptions $\Delta$ in the subarguments for minor premisses of elimination rules in the principal path.
\begin{displaymath}
  \begin{minipage}{\textwidth}
    \hbox{
      \makebox[0.3\textwidth][c]{Complementation}
      \makebox[0.4\textwidth][c]{Valid Canonical Argument}
      \makebox[0.3\textwidth][c]{Derivation}
    }
    \medskip
    \hbox{
      \makebox[0.3\textwidth][c]{
        \hbox{
          \infer*{C}{\infer{G,}{\Gamma} & \hspace{-0.7em}\Delta}
          (\deduction{c})
        }
      }

      \makebox[0.4\textwidth][c]{
        \hbox{
          \infer*{C}{\Gamma{},\Delta}
          (\deduction{v})
        }
      }

      \makebox[0.3\textwidth][c]{
        \hbox{
          \infer*{C}{\Gamma{},\Delta}
          (\deduction{d})
        }
      }
    }
  \end{minipage}
  \end{displaymath}
By \ref{def:validity}, the valid canonical argument \deduction{v} have $C$ as conclusion and, at most, $\Gamma,\Delta$ as assumptions.
We use \deduction{v} as base and assume that we already have a natural deduction derivation \deduction{d} for $C$ from $\Gamma,\Delta$ obtained by recursive application of the procedure described here to the critical subarguments of \deduction{v}.\footnote{Recall that a canonical argument is mostly already a natural deduction derivation, except for the critical subarguments (\ref{rmk:critical}).
  The termination of the recursive application of the procedure to the critical subarguments of \deduction{v} is warranted by the complexity restrictions on critical subarguments (\ref{def:validity-canonical}).
  \ref{exp:recursive} illustrates the recursive nature of the definitions.}
Then, we go through the applications of elimination rules in the principal path and construct, through a process of inversion, a natural deduction derivation of $G$ from $\Gamma$ alone.
Starting with $C$, we obtain a natural deduction derivation for each
principal sentence in the principal path until we reach $G$ (at which point we would have either discarded or discharged the auxiliary assumptions $\Delta$).
Since \ref{def:validity} ensures a valid canonical argument \deduction{v} for \emph{any} complementation \deduction{c}, we are free to consider those complementations that are more convenient for the construction of our natural deduction derivation.
We proceed case by case, where each case shows the derivation of the major premiss on the basis of derivations of the principal sentences that came before (ordered from $C$ to $G$).
For the cases of \eliminationrule{\vee} and \eliminationrule{\to}, which introduce auxiliary assumptions, I show how these assumptions can be either discarded (\eliminationrule{\vee}) or discharged (\eliminationrule{\to}).  That is, for each occurrence of \eliminationrule{\vee} and \eliminationrule{\to} in the principal path of \deduction{c}, I show how to obtain a derivation from only $\Gamma,\Delta_{*}$, where $\Delta_{*}$ stands for the auxiliary assumptions except those assumptions that are being introduced at that particular inference step.
As a result, after going through all sentences in the principal path, we obtain a derivation of $G$ which depends solely on $\Gamma$.
\begin{itemize}
\item [(\eliminationrule{\wedge})] Consider complementations by both elimination rules for conjunction
\begin{displaymath}
  \hbox{
      \makebox[0.2\textwidth][c]{
        \hbox{
          \infer{A}{\infer*{A\wedge{B}}{\infer{G,}{\Gamma}
              & \hspace{-0.7em}\Delta}}
          (\deduction{c_{1}})
        }
      }
      \makebox[0.2\textwidth][c]{
        \hbox{
          \infer{B}{\infer*{A\wedge{B}}{\infer{G,}{\Gamma}
              & \hspace{-0.7em}\Delta}}
          (\deduction{c_{2}})
        }
      }
      \makebox[0.6\textwidth][c]{
        \infer{A\wedge{B}}{
          \mbox{\hbox{
              \infer*{A}{\Gamma,\Delta}
              (\deduction{d_{1}})}}
          &
          \mbox{\hbox{
              \infer*{B}{\Gamma,\Delta}
              (\deduction{d_{2}})}}
        }
      }
  }
\end{displaymath}
From the derivations of $A$ and $B$, the conjunction $A\wedge{B}$ is derived by \introductionrule{\wedge}.

\item [(\eliminationrule{\to})] Consider a complementation
  where the minor premiss $A$ is assumed.
  \begin{displaymath}
    \hbox{
      \makebox[0.4\textwidth][c]{
        \hbox{
          \infer{B}{\infer*{A\to{B}}{\infer{G,}{\Gamma} &
              \hspace{-0.7em}\Delta_{*}} & A}
          (\deduction{c})
        }
      }
      \makebox[0.4\textwidth][c]{
        \infer{A\to{B}}{\mbox{\hbox{
              \infer*{B}{\Gamma,\Delta_{*},\discharge{A}}
              (\deduction{d})}}
        }
      }
    }
  \end{displaymath}
From a derivation of $B$ from $\Gamma,\Delta$, we apply
\introductionrule{\to} to obtain $A\to{B}$, discharging
$A$.
\item [(\eliminationrule{\bot})] Consider a complementation where $C$ is an atomic sentence which does not occur as a subsentence in either $\Gamma$ or $\Delta$.
  As a result, \eliminationrule{\bot} is the last rule applied.
  I show that \deduction{d} contain a derivation of $\bot$ which depends, at most, on $\Gamma,\Delta$.
  \begin{displaymath}
    \hbox{
      \makebox[0.4\textwidth][c]{
        \hbox{
          \infer{C}{\infer*{\bot}{\infer{G,}{\Gamma} &
              \hspace{-0.7em}\Delta}}
          (\deduction{c})
        }
      }
      \makebox[0.4\textwidth][c]{
        \hbox{
          \infer{C}{\fbox{\infer*{\bot}{\Gamma,\Delta}}}
          (\deduction{d})
        }
      }
    }
  \end{displaymath}
  By \ref{def:proper}, \deduction{d} has a principal path from one of the assumptions $\Gamma,\Delta$ to $C$.
  Since $C$ is not a subsentence of the principal assumption, it could only have been obtained by either \eliminationrule{\bot} or, possibly, a sequence of one or more applications of \eliminationrule{\vee}.\footnote{The cases where there is an application of \eliminationrule{\bot} with a complex consequence containing $C$ and then further eliminations arriving at $C$ are easily subsumed under the case where the corresponding application of \eliminationrule{\bot} has $C$ directly as a consequence.} 
  In the first case, we already have a derivation of $\bot$ from $\Gamma,\Delta$.
  In the second case, by \ref{def:canonical}, the subarguments for vertical premisses of \eliminationrule{\vee} are proper and hence have a principal path to $C$.  We work our way up applications of \eliminationrule{\vee} until we reach an application of a reductive elimination rule.
The reductive elimination rule in question can only be \eliminationrule{\bot} (\ref{lem:subsentence}).
We permute its application down the sequence of applications of \eliminationrule{\vee} and thus obtain a derivation of $\bot$ from, at most, $\Gamma{},\Delta$.
\item [(\eliminationrule{\vee})] Consider a complementation of the form below, where $A\to{C}$ and $B\to{C}$ are assumed and $C$ is atomic and does not occur as a subsentence in either $\Gamma$, $A$, $B$ or auxiliary assumptions $\Delta_{*}$ (where $\Delta_{*}$ does not contain $A\to{C}$ and $B\to{C}$).
  \begin{equation*}\tag{\deduction{c}}
    \infer{C}{\infer*{A\vee{B}}{\infer{G,}{\Gamma} & \hspace{-0.7em}\Delta_{*}} &
      \infer{C}{A\to{C} & \discharge{A}} &
      \infer{C}{B\to{C} & \discharge{B}}}
  \end{equation*}
In the derivation \deduction{d}, the conclusion $C$ could have been obtained:
\begin{enumerate}[label=(\alph*{})]
\item \label{absurdity} by \eliminationrule{\bot}.
  \begin{displaymath}
    \hbox{
      \makebox[0.4\textwidth][c]{
        \hbox{
          \infer{C}{\infer*{\bot}{\Gamma{},\Delta_{*}}}
          (\deduction{d})
        }
      }
      \makebox[0.4\textwidth][c]{
        \infer{A\vee{B}}{\infer*{\bot}{\Gamma{},\Delta_{*}}}
      }
    }
  \end{displaymath}
It is easy to derive $A\vee{B}$ instead by the same rule.
\item \label{cond} by \eliminationrule{\to} from either $A\to{C}$ or, respectively, $B\to{C}$ as principal assumption.
\begin{displaymath}
\hbox{
  \makebox[0.4\textwidth][c]{
    \hbox{
      \infer{C}{A/B\to{C} & \fbox{\infer*{A/B}{\Gamma,\Delta_{*}}}}
      (\deduction{d})
    }
  }
  \makebox[0.4\textwidth]{
    \infer{A\vee{B}}{
      \infer*{A/B}{\Gamma,\Delta_{*}}}
  }
}
\end{displaymath}
In either case, we obtain a derivation of $A\vee{B}$ from the subderivation of the horizontal minor premiss by \introductionrule{\vee}.
\item \label{disj} by a sequence of one or more applications of \eliminationrule{\vee}.
  \begin{displaymath}
    \hbox{
      \makebox[0.4\textwidth][c]{
        \hbox{
          \infer{C}{\infer*{\cdots}{\Gamma{},\Delta_{*}} & \infer{C}{A/B\to{C} & \infer*{A/B}{\Gamma{},\Delta_{*}}} &
            \infer{C}{\infer*{\bot}{\Gamma{},\Delta_{*}}}}
          (\deduction{d})
        }
      }
      \makebox[0.4\textwidth][c]{
        \infer{A\vee{B}}{\infer*{\cdots}{\Gamma{},\Delta_{*}} & \infer{A\vee{B}}{\infer*{A/B}{\Gamma{},\Delta_{*}}} &
          \infer{A\vee{B}}{\infer*{\bot}{\Gamma{},\Delta_{*}}}}
      }
    }
  \end{displaymath}
We then consider main tracks in \deduction{d} (\ref{rmk:track}) and replace each occurrence of $C$ in the segment by $A\vee{B}$ in accordance with the previous cases \ref{absurdity} and \ref{cond}.
\end{enumerate}
The resulting derivation of $A\vee{B}$ abstains from assumptions $A\to{C}$ and $B\to{C}$.
Any doubts can be dispelled by putting the derivation into normal form (\ref{thm:normal}).\qedhere
\end{itemize}
\end{proof}
\begin{theorem}[Soundness of Intuitionistic Logic]\label{thm:soundness}
  If there is a natural deduction derivation of $G$ from
  $\Gamma$, then the argument \linearargument{\Gamma}{G} is
  valid.
\end{theorem}

\begin{proof}
  Let \deduction{d} be a normal derivation of $G$ from $\Gamma$.
  Now, suppose \deduction{c} is a complementation of \deduction{d} from $\Gamma{},\Delta$ to conclusion $C$.
  I show how to obtain a valid canonical argument from $\Gamma{},\Delta$ to $C$.
  By \ref{def:proper}, $G$ is the first of a (possibly empty) sequence of major premisses of applications of elimination rules and, by \ref{thm:normal}, the last of a (possibly empty) sequence of immediate premisses of applications of introduction rules.
  By induction on the degree of $G$, we perform reductions until we obtain a deduction \deduction{v} of $C$ from $\Gamma{},\Delta$.
  By \ref{rmk:track} and induction hypotheses on its critical subarguments, \deduction{v} is a valid canonical argument for $C$ from $\Gamma{},\Delta$.
\end{proof}
\subsection{Prawitz's objection}\label{sec:prawitz-objection}
There are subtle issues involved in the treatment of critical subarguments and Dummett was not able to get his definitions completely right.
In particular, problems emerge when we consider a counterexample pointed out by \textcite[endnote~15]{prawitzPVTM}:
\begin{quote}
  ``The main fault [with Dummett's definitions] is that in a complementation of an argument, the minor premise of an implication elimination is only assumed.
  By not considering complementations where the minor premise is the end of an arbitrary argument (which is not possible to do in Dummett's definition, proceeding as it does by induction over the degree of arguments), the notion of validity becomes too weak.
  In particular, it cannot be shown that inferences by \emph{modus ponens} are valid in general, because given two valid arguments $\Pi$ and $\Sigma$ for $A\to{B}$ and $A$, respectively, there is no guarantee that the result $\Delta$ of combining them in a \emph{modus ponens} to conclude $B$ is valid.
  For an actual counterexample, we may let $B$ be atomic, $\Pi$ be simply $A\to{B}$, which is a valid argument for $A\to{B}$ from $A\to{B}$, and $\Sigma$ to be a valid argument for a nonatomic $A$ from some hypotheses of higher degree than that of $A\to{B}$.
  Then $\Delta$ is canonical argument and is its own complementation, but it is not valid ($\Sigma$ being of the same degree as $\Delta$, nor can one find another valid canonical argument for $B$ from the same hypotheses.''
\end{quote}

In his reply to Prawitz, \textcite{dummettRDP} acknowledges the problem.
There are actually two different issues brought to light by Prawitz's counterexample.
In the remainder of this section, I discuss these issues and indicate thereby the adaptations that I incorporated into the original definitions in order to avoid them.
The adaptations, although elaborate, are fully in agreement with \mbox{Dummett's} overall philosophical outlook, particularly with respect to the treatment of assumptions.  
\subsubsection{Canonical atomism}\label{sec:atomism}
\textcite[pp.~284,285]{dummettLBM} discussed an example closely related to Prawitz's counterexample.
\begin{equation}\label{eq:countexa}
  \infer{B}{A\to{B} & \infer{A}{
      (A\to{B})\to{((C\to{C})\to{A})} & A\to{B}}}
\end{equation}
Notice that \mbox{Dummett's} example is basically an instance of Prawitz's counterexample: both consist of an argument where major premiss $A\to{B}$ stands as an assumption, and where there is a subargument for minor premiss $A$ from assumptions of higher complexity than $A\to{B}$.
In Dummett's discussion, however, the minor premiss $A$ is \emph{atomic}.
This contrasts with Prawitz's counterexample where $A$ is \emph{complex}.
The difference is important because, according to Dummett's original definition (which is divided into clauses), a canonical argument, besides being proper (clause iii), must have an atomic conclusion (clause i).
Thus, the first problem revealed by Prawitz's counterexample is that, for complex $A$, there would be, in general, no canonical way to obtain $A$.

However, there is no conceptually compelling reason why canonical arguments must have atomic conclusions.
After all, we should be able to obtain also complex sentences in a canonical manner.\footnote{I suspect that Dummett only imposed the requirement of atomic conclusion on canonical arguments in order to simplify the formulation of his definition of complementation which, in general, should require that the principal path be as long as possible in order to afford a complete analysis of the conclusion.}
In order to avoid this objection, I removed the requirement of atomic conclusion from Dummett's original definition of canonical argument and adapted the definition of complementation accordingly (clause (iii) of \ref{def:complementation}).
\subsubsection{Stringency of the complexity restriction}\label{sec:stringency}
When discussing his example (\ref{eq:countexa}), \mbox{Dummett} was concerned about improper and, therefore, non-canonical subarguments for minor premisses: if these kind of subarguments could have higher complexity than the principal assumption ($A\to{B}$, in this example), the definition of validity would be in danger of circularity.
\mbox{Dummett} then presents a transformation that puts the improper subargument into proper form.
\begin{equation}\label{eq:countexb}
  \infer{B}{A\to{B} & \infer{A}{
      \infer{(C\to{C})\to{A}}{
        (A\to{B})\to{((C\to{C})\to{A})} & A\to{B}} &
      \infer{C\to{C}}{\discharge{C}}}}  
\end{equation}

Both arguments depend on the same assumptions but, in contrast with the original example (\ref{eq:countexa}), the transformed argument (\ref{eq:countexb}) displays a proper subargument for the minor premiss $A$, since there is a principal path from $(A\to{B})\to{((C\to{C})\to{A})}$ to $A$.

Apparently relying on the strength of this particular transformation alone, \mbox{Dummett} then introduces a \emph{narrow} notion of validity for canonical arguments which restricts improper arguments for minor premisses to those of \emph{strictly lower} degree, where the degree of an argument is the highest among the degrees of its assumptions and conclusion.

Although the transformation worked for that particular example, it is inadequate in general, at least if \mbox{Dummett's} notion of degree of an argument is used as complexity measure.
Consider, for instance, the following proper argument:
\begin{displaymath}
  \infer{\bot}{
    \lnot{(A\vee{\lnot{A}})} &
      \infer{A\vee{\lnot{A}}}{\lnot{(A\vee{\lnot{A}})}}
  }
\end{displaymath}
The degree of the minor subargument is \emph{equal} to the degree of the principal assumption $\lnot{(A\vee{\lnot{A}})}$.
In fact, $\lnot{(A\vee{\lnot{A}})}$ occurs again as an assumption in the minor subargument.
The fact that the minor subargument cannot be put into a proper form becomes clear when we replace it by its normal derivation in intuitionistic logic:
\begin{displaymath}
  \infer{\bot}{
    \lnot{(A\vee{\lnot{A}})} &
    \infer{A\vee{\lnot{A}}}{
      \infer[\text{\scriptsize{(1)}}]{\lnot{A}}{
        \infer{\bot}{
          \lnot{(A\vee{\lnot{A}})} &
          \infer{A\vee{\lnot{A}}}{A^{(1)}}
        }
      }
    }
  }
  \end{displaymath}
\mbox{Dummett's} complexity restriction, as originally formulated, is therefore too stringent.
The approach suggested by \textcite{prawitzPVTM} and \textcite{schroederheisterPTVBER} avoids this problem by dealing primarily with closed proofs, where the conclusion provides the adequate complexity measure, since there are no open assumptions.
Their approach thus differs unequivocally from \mbox{Dummett's}, especially with respect to the treatment of assumptions.
Dummett's core approach was maintained through the adoption of an adequate complexity measure (\ref{def:complexity}), one that preserves soundness (\ref{thm:complexity}), instead of his original notion of degree of an argument.

\section{A decision procedure}\label{sec:decision-procedure}
In order to illustrate the definitions and give some intuition about the construction described in the proof of \ref{thm:completeness}, it can be useful to work though some examples.
The examples are meant to be an overall intuitive illustration of how the proof-theoretic definitions evaluate the validity of arguments.
They are presented in the framework of a decision procedure that can be read off from the definitions.

The idea behind \ref{def:validity} is roughly that an argument is valid if, whatever we can obtain canonically from the conclusion, could as well be obtained from the assumptions.
A procedure to evaluate validity can therefore be divided into two
parts:
\begin{description}
\item The \texttt{complementation} process determines what can be obtained from
  the conclusion.
\item The \texttt{search} process looks for a way to obtain the same thing from
  the assumptions.
\end{description}
Both \texttt{complementation} and \texttt{search} can employ only elimination rules --- there are no introduction rules available.
In line with \ref{def:proper}, they are based on a similar method (let us call it \emph{analysis}) of applying elimination rules to a sentence, taken as major premiss, until an atomic sentence is obtained (clause \ref{cond-atomic} of \ref{def:complementation}).
Thus, in the \texttt{complementation} process, the conclusion of the argument is analysed, in order to see what atomic conclusions can be obtained (possibly under some additional auxiliary assumptions).
In the \texttt{search} process, the assumptions are then analysed (one by one), in order to evaluate whether the same atomic conclusions can be obtained.

In the \texttt{complementation} process, the following simplifications are adopted, without loss of generality, with respect to $\vee$E and $\to$E (in agreement with the corresponding cases in the proof of \ref{thm:completeness}):
\begin{itemize}
\item [($\to$E)] the minor premiss is assumed.
  \begin{displaymath}
    \infer{B}{A\to{B} & A}
  \end{displaymath}
  Here, $A$ is an additional assumption and will be available to
  \texttt{search}.
\item [($\vee$E)] applications are ``flattened'' with the help of implication.
  \begin{displaymath}
    \infer{C}{A\mathbin{\vee}{B} &
      \infer{C}{A\mathbin{\to}{C} & [A]} &
      \infer{C}{B\mathbin{\to}{C} & [B]}}
  \end{displaymath}
  In order to maintain generality, $C$ stands for a sentence that does not occur as a subsentence either in the assumptions or the conclusion.
  Here, $A\mathbin{\to}{C}$ and $B\mathbin{\to}{C}$ are assumed and will be available to \texttt{search}.
\item [($\bot$E)] applications are abstained.
The \texttt{search} will then target $\bot$.  Notice that these simplifications are limited to the \texttt{complementation} process and do \emph{not} carry over to the \texttt{search} process where, naturally, applications of $\bot$E are not abstained.
\end{itemize}
\begin{example}\label{exp:simple}
A definition of validity is expected to provide precise criteria for the validity of arguments and, for our definitions in particular, these criteria are supposed to resort to elimination rules only (without assistance from introduction rules).
Consider a simple, but not trivial, argument
\begin{displaymath}
  \infer{(A\mathbin{\to}B)\mathbin{\wedge}(A\mathbin{\to}C)}{
    A\mathbin{\to}(B\mathbin{\wedge}C)}
\end{displaymath}
and let us evaluate its validity with respect to our definitions.

First, we investigate what can be obtained canonically from the conclusion by means of \texttt{complementation}:
\begin{displaymath}
  \begin{array}{p{0.6\textwidth}p{0.3\textwidth}}
    \infer{B}{
      \infer{A\mathbin{\to}B}{
        \infer{(A\mathbin{\to}B)\mathbin{\wedge}(A\mathbin{\to}C)}{
          A\mathbin{\to}(B\mathbin{\wedge}C)}
      } & A}
    &
    \infer{C}{
      \infer{A\mathbin{\to}C}{
        \infer{(A\mathbin{\to}B)\mathbin{\wedge}(A\mathbin{\to}C)}{
          A\mathbin{\to}(B\mathbin{\wedge}C)}
      } & A}
  \end{array}
\end{displaymath}
There are two complementations, with conclusions $B$ and $C$, respectively, and the assumptions $A\mathbin{\to}(B\mathbin{\wedge}C)$ and $A$.
In order to establish validity, we must now find canonical arguments from $A\mathbin{\to}(B\mathbin{\wedge}C)$ and $A$ to $B$, and from $A\mathbin{\to}(B\mathbin{\wedge}C)$ and $A$ to $C$.\footnote{The assumptions of the complementations happen to be the same in this example.
  In the general case, however, they have to be considered separately, e.g.\,{} each complementation has their own assumptions and conclusion.
  In order to establish validity, we must then show that the conclusion of the complementation can be obtained from the assumptions of the complementation, \emph{for every complementation}.}
The \texttt{search} for these canonical arguments can be done mechanically by analysing the assumptions one by one, where some heuristics could be employed to sort out the most promising candidates.
In this example, we have few assumptions and don't need much heuristics to see that $A\mathbin{\to}(B\mathbin{\wedge}C)$ is the best candidate:
\begin{displaymath}
  \begin{array}{p{0.6\textwidth}p{0.3\textwidth}}
    \infer{B}{
      \infer{B\mathbin{\wedge}C}{
        A\mathbin{\to}(B\mathbin{\wedge}C) & A}}
    &
    \infer{C}{
      \infer{B\mathbin{\wedge}C}{
        A\mathbin{\to}(B\mathbin{\wedge}C) & A}}
  \end{array}
\end{displaymath}
The procedure is thus revealed to be strong enough to validate, not only the introduction rules on the basis of the elimination rules, but also complex arguments whose derivation would require both eliminations \emph{and} introduction rules.
\end{example}
\begin{example}\label{exp:recursive}
Regarded as a decision algorithm, the procedure for evaluation of validity based on elimination rules is not so straightforward and uncomplicated as \ref{exp:simple} makes it out to be.
In the general case, the procedure may involve recursion and backtracking.
The \texttt{search} process can deliver candidates with critical subarguments, which would demand a recursive call to evaluate their legitimacy (\ref{def:validity-canonical}).
If unsuccessful, the process backtracks and tries out the analysis on a different assumption.
Consider, for instance, the argument
\begin{displaymath}
  \infer{\lnot\lnot{}(\lnot{A}\mathbin{\vee}B)}{
    A\mathbin{\to}\lnot\lnot{B}}
\end{displaymath}
The \texttt{complementation} below stops at the conclusion $\bot$, before an application of \eliminationrule{\bot}, in accordance with the aforementioned simplifications to the complementation process.
\begin{equation*}\tag{C1}
  \infer{\bot}{
    \infer{\lnot\lnot{}(\lnot{A}\mathbin{\vee}B)}{
      A\mathbin{\to}\lnot\lnot{B}} &
    \lnot{}(\lnot{A}\mathbin{\vee}B)}
\end{equation*}
The \texttt{search} process has assumptions $A\to{\lnot\lnot{B}}$ and $\lnot{(\lnot{A}\vee{B})}$ to try out in order to obtain $\bot$.
For simplicity of exposition, we heuristically select $\lnot{(\lnot{A}\vee{B})}$, but could as well have unsuccessfully tried $A\to{\lnot\lnot{B}}$ out and backtracked here.
\begin{displaymath}
  \infer{\bot}{
    \lnot{}(\lnot{A}\mathbin{\vee}B) & \lnot{A}\mathbin{\vee}B}
\end{displaymath}
Now, notice that $\lnot{A}\mathbin{\vee}B$ itself is not available among our assumptions.
Therefore we \emph{presume} that $\lnot{A}\mathbin{\vee}B$ can, in fact, be obtained from the assumptions that \emph{are} available to us, and recall the procedure recursively on the critical subargument enclosed in a box below.\footnote{Notice that $\lnot{}(\lnot{A}\mathbin{\vee}B)$ appears twice: as major premiss and also as an assumption of the critical subargument.
This cannot be avoided in general and is related to the problem of contraction in the search for proofs in the sequent calculus \autocite{dyckhoffCFSCIL}.}
\begin{equation*}\tag{S1}
  \infer{\bot}{
    \lnot{}(\lnot{A}\mathbin{\vee}B) &
    \fbox{
      \infer{\lnot{A}\mathbin{\vee}B}{
        A\mathbin{\to}\lnot\lnot{B} &
        \lnot{}(\lnot{A}\mathbin{\vee}B)}}}
\end{equation*}
In the \texttt{complementation} of our recursive call, we again adhere to aforementioned simplifications and use $C$ as the conclusion of $\vee$E since it does not occur anywhere else.
\begin{equation*}\tag{C2}
  \infer{C}{
    \infer{\lnot{A}\mathbin{\vee}B}{
      A\mathbin{\to}\lnot\lnot{B} &
      \lnot{}(\lnot{A}\mathbin{\vee}B)} &
    \infer{C}{\lnot{A}\mathbin{\to}C & [\lnot{A}]} &
    \infer{C}{B\mathbin{\to}C & [B]}}
\end{equation*}
In order to obtain our foreign $C$, the \texttt{search} must either (1) obtain $\bot$, and thereby $C$, or (2) obtain one of the disjuncts and thereby obtain $C$ from the corresponding assumption, either $\lnot{A}\mathbin{\to}C$ or $B\mathbin{\to}C$, or yet (3) obtain $C$ by \eliminationrule{\vee} from a disjunctive principal sentence, whereby we may use the disjuncts as additional assumptions on the \texttt{search} for proper subarguments for the respective vertical minor premisses.
We examine the second option and choose assumption $\lnot{A}\mathbin{\to}C$.
The other one may be discarded.
\begin{equation*}\tag{S2}
  \infer{C}{
    \lnot{A}\mathbin{\to}C &
    \fbox{
      \infer{\lnot{A}}{
        A\mathbin{\to}\lnot\lnot{B} &
        \lnot{}(\lnot{A}\mathbin{\vee}B)}}}
\end{equation*}
Our next recursive step reveals an important aspect of the definitions.
Consider the next \texttt{complementation}.
\begin{equation*}\tag{C3}
  \infer{\bot}{
    \infer{\lnot{A}}{
      A\mathbin{\to}\lnot\lnot{B} &
      \lnot{}(\lnot{A}\mathbin{\vee}B)} & A}
\end{equation*}
In the candidate (S3) below, if we were to retain all the assumptions available for the next recursive call, that is, if $A\to{\lnot\lnot{B}}$ and $A$ where \emph{both} passed as assumptions to the critical subargument enclosed in a box, we would be in danger of running into a vicious circle (a loop):
after the \texttt{complementation} (C4) below, the candidate (S2) above could have been considered by the \texttt{search}.  Indeed, by \ref{def:complexity}, the argument from $\{A\to{\lnot\lnot{B}}, A, \lnot{(\lnot{A}\vee{B})}\}$ to $\lnot{B}$ has higher complexity than the critical subargument in (S2), because its conclusion $\lnot{B}$ has the same degree than $\lnot{A}$ and it has $A$ as an additional assumption.  Therefore, for the particular case with $A\to{\lnot\lnot{B}}$ as principal assumption, the \texttt{search} must consider only candidates where $A$ or some other assumptions is left out of the critical subargument, on pain of violating the complexity restriction.
As it turns out, we do not need $A\to{\lnot\lnot{B}}$ either.
\begin{equation*}\tag{S3}
  \infer{\bot}{
    \infer{\lnot\lnot{B}}{
      A\mathbin{\to}\lnot\lnot{B} & A} &
    \fbox{
      \infer{\lnot{B}}{\lnot{}(\lnot{A}\mathbin{\vee}B)}}}
\end{equation*}
More recursion.
\begin{equation*}\tag{C4}
  \infer{\bot}{
    \infer{\lnot{B}}{\lnot{}(\lnot{A}\mathbin{\vee}B)} & B}
\end{equation*}
We see the complexity restriction at work again in the candidate offered by the \texttt{search} below (notice that $\lnot{}(\lnot{A}\mathbin{\vee}B)$ is left out of the critical subargument).
\begin{equation*}\tag{S4}
  \infer{\bot}{
    \lnot{}(\lnot{A}\mathbin{\vee}B) &
    \fbox{
      \infer{\lnot{A}\mathbin{\vee}B}{B}}}
\end{equation*}
I think that the procedure should be clear enough by now for us to omit the last recursive call.
\end{example}
The construction described in \ref{thm:completeness} can be applied to the canonical arguments produced by \texttt{complementation} and \texttt{search} in order to obtain a derivation.
\begin{displaymath}
  \infer[\text{\scriptsize{(1),(3)}}]{
    \lnot\lnot{(\lnot{A}\vee{B})}}{
    \infer{\bot}{
      \lnot{(\lnot{A}\vee{B})}^{(1)} &
      \infer{\lnot{A}\vee{B}}{
        \infer[\text{\scriptsize{(2)}}]{\lnot{A}}{
          \infer{\bot}{
            \infer{\lnot\lnot{B}}{
              A\to{\lnot\lnot{B}} &
              A^{(2)}} &
            \infer[\text{\scriptsize{(4)}}]{\lnot{B}}{
              \infer{\bot}{
                \lnot{(\lnot{A}\vee{B})}^{(3)} &
                \infer{\lnot{A}\vee{B}}{
                  B^{(4)}}
              }
            }
          }
        }
      }
    }
  }
\end{displaymath}
The derivation contains four tracks.  If we order the tracks from one to four and divide them into their analytic and synthetic parts, they correspond roughly to the complementation and search processes of the procedure: (C1), [(C2), (C3)], (C4) and (C5) (omitted) correspond to the synthetic parts of tracks 1, 2, 3 and 4; (S1), (S3) and (S4) correspond to the analytic parts of tracks 1, 2 and 3 (the analytic part of track 4 is empty).
The simplifications adopted with respect to \eliminationrule{\vee} in the \texttt{complementation} process resulted in a dedicated recursive step for applications of \introductionrule{\vee} in the derivation (in track 2, (C2) and (S2); in track 4, (C5) and (S5)).
This seems a reasonable exchange against the achieved separation between the processes and deterministic character of the \texttt{complementation}.
\section{Discussion}
Dummett's pragmatist justification procedure rejects a widely accepted dog\-ma of proof-theoretic semantics: \emph{the primacy of the categorical over the hypothetical} or, as it is also called, \emph{the placeholder view of assumptions} \autocite{schroederheisterPTVMTC,schroederheisterCH}.
According to this view, assumptions are placeholders for closed proofs and thus hypothetical reasoning (reasoning from assumptions) are explained in terms of categorical reasoning (proofs without assumptions).
In contrast, the pragmatist proof-theoretic notion of canonical argument considers arguments from assumptions as primary and not to be explained away in terms of closed arguments (proofs).

The prevalence of the placeholder view of assumptions has imparted greater emphasis to assertions in detriment of other speech acts like denial or supposition.
In contrast, the approach presented in this paper indicate ways to widen the conceptual arena and thus accommodate interesting alternative points of view.

\paragraph{Acknowledgements.} 
I received very helpful comments and suggestions from members of the research group ``Logik und Sprachtheorie'' in T\"{u}bingen.
I am indebted in particular to Peter Schroeder-Heister, who suggested to me the idea of pragmatist canonical arguments with complex conclusions, and to Luca Tranchini, who read and made comments on early drafts.
This work was supported by the Deutscher Akademischer Austauschdienst (DAAD), grant number 91562976.


\begin{thebibliography}{22}
\providecommand{\natexlab}[1]{#1}
\providecommand{\url}[1]{\texttt{#1}}
\expandafter\ifx\csname urlstyle\endcsname\relax
  \providecommand{\doi}[1]{doi: #1}\else
  \providecommand{\doi}{doi: \begingroup \urlstyle{rm}\Url}\fi

\bibitem[Dershowitz and Manna(1979)]{dershowitzPTMO}
Nachum Dershowitz and Zohar Manna.
\newblock Proving termination with multiset orderings.
\newblock \emph{Communications of the ACM}, 22\penalty0 (8):\penalty0 465--476,
  1979.

\bibitem[Dummett(1975)]{dummettPBIL}
Michael Dummett.
\newblock The philosophical basis of intuitionistic logic.
\newblock In H.E. Rose and J.C. Shepherdson, editors, \emph{{L}ogic
  {C}olloquium '73}, volume~80 of \emph{Studies in Logic and the Foundations of
  Mathematics}, pages 5--40. North Holland, Amsterdam, 1975.

\bibitem[Dummett(1991)]{dummettLBM}
Michael Dummett.
\newblock \emph{{T}he {L}ogical {B}asis of {M}etaphysics}.
\newblock Harvard University Press, Cambridge, Massachusetts, 1991.

\bibitem[Dummett(2007)]{dummettRDP}
Michael Dummett.
\newblock Reply to {D}ag {P}rawitz.
\newblock In Randall~E. Auxier and Lewis~Edwin Hahn, editors, \emph{The
  Philosophy of Michael Dummett}, volume XXXI of \emph{The Library of Living
  Philosophers}, chapter~13, pages 482--489. Open Court Publishing Company,
  Chicago and La Salle, Illinois, 2007.

\bibitem[Dyckhoff(1992)]{dyckhoffCFSCIL}
Roy Dyckhoff.
\newblock Contraction-free sequent calculi for intuitionistic logic.
\newblock \emph{The Journal of Symbolic Logic}, 57\penalty0 (3):\penalty0
  795--807, 9 1992.

\bibitem[Gentzen(1934)]{gentzenULS}
Gerhard Gentzen.
\newblock Untersuchungen {\"{u}}ber das logische schlie{\ss}en {I}.
\newblock \emph{Mathematische Zeitschrift}, 39\penalty0 (2):\penalty0 176--210,
  1934.

\bibitem[Goldfarb(2016)]{goldfarbDPTJLL}
Warren Goldfarb.
\newblock On {D}ummett's proof-theoretic justification of logical laws.
\newblock In Thomas Piecha and Peter Schroeder-Heister, editors,
  \emph{{A}dvances in {P}roof-{T}heoretic {S}emantics}, volume~43 of
  \emph{Trends in Logic}, pages 195--210. Springer, 2016.

\bibitem[Piecha(2016)]{piechaCPTS}
Thomas Piecha.
\newblock Completeness in proof-theoretic semantics.
\newblock In Thomas Piecha and Peter Schroeder-Heister, editors, \emph{Advances
  in Proof-Theoretic Semantics}, volume~43 of \emph{Trends in Logic}, pages
  231--251. Springer, 2016.

\bibitem[Piecha et~al.(2015)Piecha, Sanz, and
  Schroeder-Heister]{piechaetalFCPTS}
Thomas Piecha, Wagner~{de Campos} Sanz, and Peter Schroeder-Heister.
\newblock Failure of completeness in proof-theoretic semantics.
\newblock \emph{Journal of Philosophical Logic}, 44\penalty0 (3):\penalty0
  321--335, 2015.

\bibitem[Prawitz(1965)]{prawitzND}
Dag Prawitz.
\newblock \emph{{N}atural {D}eduction: {A} {P}roof-{T}heoretical {S}tudy}.
\newblock Almqvist \& Wiksell, Stockholm, 1965.

\bibitem[Prawitz(1971)]{prawitzIRPT}
Dag Prawitz.
\newblock Ideas and results in proof theory.
\newblock In \emph{{P}roceedings of the {S}econd {S}candinavian {L}ogic
  {S}ymposium}, volume~63 of \emph{Studies in Logic and the Foundations of
  Mathematics}, pages 235--307, 1971.

\bibitem[Prawitz(1973)]{prawitzTFGPT}
Dag Prawitz.
\newblock Towards a foundation of a general proof theory.
\newblock In Patrick Suppes, Leon Henkin, Athanase Joja, and G.~C. Moisil,
  editors, \emph{{L}ogic, {M}ethodology and {P}hilosophy of {S}cience {IV}},
  volume~74 of \emph{Studies in Logic and the Foundations of Mathematics},
  pages 225--250, Amsterdam, 1973. North Holland.

\bibitem[Prawitz(1998)]{prawitzTOFVPV}
Dag Prawitz.
\newblock Truth and objectivity from a verificationist point of view.
\newblock In H.~G. Dales and G.~Oliveri, editors, \emph{Truth in Mathematics},
  pages 41--51. Claredon Press, Oxford, 1998.

\bibitem[Prawitz(2006)]{prawitzMAP}
Dag Prawitz.
\newblock Meaning approached via proofs.
\newblock \emph{Synthese}, 148\penalty0 (3):\penalty0 507--524, 2006.

\bibitem[Prawitz(2007)]{prawitzPVTM}
Dag Prawitz.
\newblock Pragmatist and verificationist theories of meaning.
\newblock In Randall~E. Auxier and Lewis~Edwin Hahn, editors, \emph{{T}he
  {P}hilosophy of {M}ichael {D}ummett}, volume XXXI of \emph{The Library of
  Living Philosophers}, chapter~13, pages 455--481. Open Court Publishing
  Company, Chicago and La Salle, Illinois, 2007.

\bibitem[Prawitz(2014)]{prawitzAGPTCCILR}
Dag Prawitz.
\newblock An approach to general proof theory and a conjecture of a kind of
  completeness of intuitionistic logic revisited.
\newblock In Luiz~Carlos Pereira, Edward~Hermann Haeusler, and Val\'{e}ria {de
  Paiva}, editors, \emph{{A}dvances in {N}atural {D}eduction}, volume~39 of
  \emph{Trends in Logic}, pages 269--279. Springer, Dordrecht, 2014.

\bibitem[Sandqvist(2009)]{sandqvistCLWB}
Tor Sandqvist.
\newblock Classical logic without bivalence.
\newblock \emph{{A}nalysis}, 69\penalty0 (2):\penalty0 211--218, 4 2009.

\bibitem[Sanz et~al.(2014)Sanz, Piecha, and Schroeder-Heister]{sanzetalCSARVPL}
Wagner~{de Campos} Sanz, Thomas Piecha, and Peter Schroeder-Heister.
\newblock Constructive semantics, admissibility of rules and the validity of
  {P}eirce's law.
\newblock \emph{Logic Journal of the IGPL}, 22\penalty0 (2):\penalty0 297--308,
  2014.

\bibitem[Schroeder-Heister(2006)]{schroederheisterVCPTS}
Peter Schroeder-Heister.
\newblock Validity concepts in proof-theoretic semantics.
\newblock \emph{Synthese}, 148\penalty0 (3):\penalty0 525--571, 2006.

\bibitem[Schroeder-Heister(2008)]{schroederheisterPTVMTC}
Peter Schroeder-Heister.
\newblock Proof-theoretic versus model-theoretic consequence.
\newblock In Michal Peli\v{s}, editor, \emph{{T}he {L}ogica {Y}earbook 2007},
  pages 187--200. Filosofia, Prague, 2008.

\bibitem[Schroeder-Heister(2012)]{schroederheisterCH}
Peter Schroeder-Heister.
\newblock The categorical and the hypothetical: A critique of some fundamental
  assumptions of standard semantics.
\newblock \emph{Synthese}, 187\penalty0 (3):\penalty0 925--942, 2012.

\bibitem[Schroeder-Heister(2015)]{schroederheisterPTVBER}
Peter Schroeder-Heister.
\newblock Proof-theoretic validity based on elimination rules.
\newblock In Edward~Hermann Hauesler, Wagner de~Campos~Sanz, and Bruno Lopes,
  editors, \emph{Why is this a Proof? (Festschrift for Luiz Carlos Pereira)},
  volume~27 of \emph{Tributes}, pages 159--176. College Publications, London,
  UK, 2015.

\end{thebibliography}

\end{document}